\documentclass[12pt]{amsart}
\newtheorem{theorem}{Theorem}[section]
\newtheorem{lemma}[theorem]{Lemma}

 \theoremstyle{definition}
\newtheorem{definition}[theorem]{Definition}

\newtheorem{claim}[theorem]{Claim}
\theoremstyle{remark}
\newtheorem{remark}[theorem]{Remark}

\numberwithin{equation}{section}

\begin{document}
\setlength{\baselineskip}{1.2\baselineskip}

\title[A Pogorelov estimate and a Liouville type theorem to parabolic $k$-Hessian equations]
{A Pogorelov estimate and a Liouville type theorem to parabolic $k$-Hessian equations}

\author{Yan He, Haoyang Sheng, Ni Xiang}
\address{Faculty of Mathematics and Statistics, Hubei Key Laboratory of Applied Mathematics, Hubei University,  Wuhan 430062, P.R. China}
\email{helenaig@hotmail.com; 907026694@qq.com; nixiang@hubu.edu.cn}
\thanks{This research was supported by funds from Hubei Provincial Department of Education Key Projects D20171004.}

\begin{abstract}
We consider Pogorelov type estimates and
Liouville type theorems
to parabolic $k$-Hessian equations of the form
$-u_t \sigma_k (D^2u) =1$
in $\mathbb{R}^n\times (-\infty, 0]$.
We derive that any \textbf{$k+1$-convex-monotone}
solution to
$-u_t \sigma_k  (D^2u) =1$
 when $u(x,0)$ satisfies a quadratic growth and $0<m_1\le -u_t\le m_2$
must be a linear function of $t$ plus
a quadratic polynomial of $x$.


\end{abstract}

\maketitle
\bigskip
\section{Introduction}

In this paper, we derive a Liouville type theorem for parabolic k-Hessian equations
\begin{equation}\label{4.10}
-u_t \sigma_k (D^2u) =1, \ \textrm{in}\  \mathbb{R}^n\times (-\infty, 0].
\end{equation}
Namely, any \textbf{$k+1$-convex-monotone}
solution of \eqref{4.10}
 with a quadratic growth and $0<m_1\le -u_t\le m_2$,
must be a linear function of $t$ plus
a quadratic polynomial of $x$.

To obtain the Liouville type theorem, the key points are Pogorelov estimates in our method.
Thus, we consider
the following
equations
\begin{equation}\label{1.100}
\left\{
\begin{array}{rl}
-u_t \sigma_k (D^2u) =1, \ \textrm{in}\  \Omega,\\
u=0, \ \textrm{on}\  \partial \Omega,
\end{array}
\right.
\end{equation}
where $D^2u\in \Gamma^{k+1}$,
$0<m_1\le- u_t\le m_2$.
Here  $\Omega\subset \mathbb{R}^n\times(-\infty,0]$
is a bounded domain and $t\leq 0,$
$\Omega(t) =\{x\in \mathbb{R}^n| (x,t)\in \Omega\},$
 $t_0=\inf\{t\leq0|\Omega(t)\neq\emptyset\}$. The parabolic boundary $\partial\Omega$
is defined by
$$\partial\Omega=(\overline{\Omega(t_0)}\times{t_0})\cup \bigcup_{t\leq 0}(\partial\Omega(t)\times \{t\}),$$
where $\overline{\Omega(t_0)}$ denotes the closure of $\Omega(t_0)$ and $\partial \Omega(t)$ denotes the boundary of $\Omega(t)$.
The k-th elementary symmetric polynomial
is denoted by $\sigma_k$:
\[\sigma_k(\lambda)=\sum_{1\le i_1<\cdots<i_k\le n}\lambda_{i_1}\cdots\lambda_{i_k}.\]
$\sigma_k(D^2u)$ means $\sigma_k$
is applied to the eigenvalues of
$D^2u$.
Let $\Gamma^k$ be an open convex cone in $\mathbb{R}^n$:
\[\Gamma^k=\{\lambda=(\lambda_1,\cdots,\lambda_n)\in\mathbb{R}^n|
\sigma_1(\lambda)>0,\cdots,\sigma_k(\lambda)>0 \}.\]
Here the function $u = u(x, t) : \mathbb{R}^n \times (\infty , 0] \rightarrow \mathbb{R} $ is said to be $k$-convex
if the eigenvalues of $D^2u$ lie
in $\Gamma^k$. Moreover, it is said to be
\textbf{$k$-convex-monotone} if it is $k$ convex in $x $ and non-increasing in $t$.
The quadratic growth means that there are
some positive constants $b, c$ and sufficiently large $R$, such that,
\begin{eqnarray}\label{gro}
u(x)\geq b|x|^2-c, \quad \mbox{for} \quad |x|\geq R.
\end{eqnarray}

A priori estimates for elliptic $k$-Hessian
equations
\begin{equation}\label{070202}\sigma_k(D^2u)=f\end{equation} have been studied intensively by many authors.
In Chou-Wang \cite{CW}, the authors
got interior gradient and second order estimates when $f$ depends on $x,\ u$.
 Warren-Yuan \cite{WY}  obtained $C^2$ interior estimates in the case of equations
  $\sigma_2(D^2u) = 1$ in $\mathbb{R}^3$, which originated from special Lagrangian geometry. Guan-Qiu \cite{GQ19} established interior $C^2$ estimates for solutions of the prescribing scalar curvature
equations and 2-Hessian equations under additional assumption that $\sigma_3(D^2 u)> -A$ for some constant $A>0.$
The purely interior $C^2$ estimates for semi-convex solutions of above equation have been obtained by McGonagle-Song-Yuan \cite{MSY} recently.
For $k\geq 2$, Li-Ren-Wang \cite{LRW}
established Pogorelov estimates
under the condition $k+1$-convex, when $f$ depends on $x, u, Du$.

Our paper is based on the work of
Li-Ren-Wang  \cite{LRW}. Firstly, We extend the Pogorelov estimate
from elliptic Hessian equations
to parabolic Hessian equations.
We have obtained the following Pogorelov type estimates.

\begin{theorem}\label{estimate7}
Let $u$ be a \textbf{$k+1$-convex-monotone} solution
of (\ref{1.100}) satisfying
$0<m_1\le -u_t\le m_2$.
Then there exists a
positive constant
$\beta$ sufficiently large such that
\begin{equation}(-u)^{\beta}\Delta u\le C,\label{1.5}\end{equation}
$C$ depends on the diameter of
$\Omega(t)$, $m_1$, $m_2$, $k$ and $\sup|u|$.
\end{theorem}

For $k=2$, we can decrease the power in (\ref{1.5}) and
improve the estimates
as follows.

\begin{theorem}\label{estimate1}
Let $u$ be a $3$-convex-monotone solution
of the following equation (\ref{1.101}) satisfying
$0<m_1\le -u_t\le m_2$.

\begin{equation}\label{1.101}
\left\{
\begin{array}{rl}
-u_t \sigma_2  (D^2u) =1, \ \textrm{in}\  \Omega,\\
u=0, \ \textrm{on}\  \partial \Omega.
\end{array}
\right.
\end{equation}
Then
\begin{equation}(-u)^8\Delta u\le C,\label{100}\end{equation}
$C$ depends on the diameter of
$\Omega(t)$, $m_1$, $m_2$ and $\sup|u|$.
\end{theorem}

These type of interior estimates are important for
existence of isometric embedding of non-compact surfaces and for Liouville type theorems.
There has been much activities on Liouville type theorems for elliptic $k$-Hessian equations.
In 2003, Bao-Chen-Guan-Ji \cite{BCGJ03}~studied the
Liouville theorem to ~
\begin{equation}\frac{\sigma_k(D^2u)}{\sigma_l(D^2u)}=1,\ (k>l).\label{eqkl}\end{equation}
They proved that entire convex solutions
of the equation (\ref{eqkl}) with a quadratic
growth are quadratic polynomials.
 In 2010, Chang-Yuan \cite{CY10}~
 considered
 \begin{equation}\sigma_2(D^2u)=1.\label{CYEQ}\end{equation} and obtained that
 the entire solution to (\ref{CYEQ})
 is quadratic polynomial
 if
~$$D^2u\ge[\delta-\sqrt{\frac{2}{n(n-1)}}]I,$$
where $\delta>0$.
In 2016,
Li-Ren-Wang \cite{LRW}~considered
 ~$\sigma_k(D^2u)=1$ for general $k$.
They obtained that global \textbf{$k+1$-convex}
solutions
with a quadratic growth are quadratic
polynomials. Chen-Xiang \cite{CX} improved the condition
from
\textbf{$k+1$-convex} to \textbf{$k$-convex} for
$k=2$ under $\sigma_3(D^2u)\geq -A$. Especially, for $n=3$, $\sigma_3(D^2 u)\geq -A$ can be redundant. Then He-Sheng-Xiang \cite{HSX} removed the condition
 $\sigma_3(D^2u)\geq -A$
for $2$-Hessian equations in general dimension $n$.

However, as far as we know, Liouville type
theorems for parabolic fully nonlinear equations are known most for parabolic Monge-Amp\`ere equations.
Guti$\acute{e}$rrez-Huang  \cite{GH} extended Theorem
of J$\ddot{o}$rgens, Calabi, and Pogorelov to parabolic Monge-Amp\`ere equations.
Xiong-Bao  \cite{XB11} obtained Liouville theorems for $$u_t=(\det D^2u)^{1/n}.$$ Zhang-Bao-Wang  \cite{ZBW} extend the theorem of Caffarelli and Li  \cite{CL} to parabolic Monge-Amp\`ere equation $$-u_t\det D^2u =f ,$$ and obtain asymptotic behavior at infinity. And along the line of approach in their paper, other parabolic Monge-Amp\`ere equations can be also treated.
 For general $k$, Nakamori S. and Takimoto K.\cite{NT} studied the bernstein type theorem for parabolic $k$-Hessian equations when the entire solution $u$ was convex-monotone.
Recently, He-Pan-Xiang \cite{HPX} prove that the \textbf{$2$-convex-monotone}
solutions
with $\sigma_3(D^2u)>-A$, $m_1\le -u_t\le m_2$ and a quadratic growth
must be a linear function of $t$ plus
a quadratic polynomial of $x$ when $k=2$.

Then using Theorem \ref{estimate7}, we have established the following Liouville type theorem for parabolic $k$-Hessian equations.

\begin{theorem}\label{maintheorem}
Let $u$ be a \textbf{$k+1$-convex-monotone} solution
of (\ref{4.10}), $u(x,0)$ satisfying a quadratic growth, and
$0<m_1\leq -u_t\leq m_2$. Then $u$ has the form
$u(x, t)=-mt+ p(x)$
where $m>0$ and $p(x)$ is a quadratic
polynomial.
\end{theorem}

This paper is organized as follows.
We start with some notations
and Lemmas in section 2.
In section 3 we prove a Pogorelov estimate for the \textbf{$k+1$-convex-monotone} solutions to parabolic $k$-Hessian equation (\ref{1.100}).
A Pogorelov estimate for the \textbf{2-convex-monotone} solutions to parabolic 2-Hessian equation (\ref{1.101}) is given in section 4.
The proof of Liouville Theorem (Theorem\ref{maintheorem})
is given in section 5.

\section{Preliminaries}
Throughout this paper, we use the Einstein summation
convention and denote by $\lambda=(\lambda_1,\cdots,\lambda_n)$
the eigenvalues of $D^2u$.
To begin this section, we introduce some notations.

\begin{definition}\label{def2}Let $\lambda=(\lambda_1,\cdots,\lambda_n)\in \mathbb{R}^n$.

(1)\[\sigma_l(\lambda|i)=\sigma_{l}(\lambda)\big|_{\lambda_i=0}.\]
 ~$\sigma_{k-1}(\lambda|i)$~is also denoted by~$\sigma_k^{ii}$.

(2)\[\sigma_l(\lambda|pq)=\sigma_l(\lambda)
\big|_{\lambda_p=\lambda_q=0}.\]
 ~$\sigma_{k-2}(\lambda|pq)$~is also denoted by~$\sigma_k^{pp,qq}$.
\end{definition}

The following Lemmas will be used in the proof for Pogorelov estimates.

\begin{lemma}(See  \cite{S})\label{lemma2}
 Suppose $\lambda\in \Gamma^k$.
For $0\le l<k\le n, 0\le s<r\le n, k\ge r, l\ge s$, the following
is the generalized Newton-MacLaurin inequality
\[\big[\frac{\sigma_k(\lambda)/C^k_n}{\sigma_l(\lambda)/C^l_n}\big]^{\frac{1}{k-l}}
\le \big[\frac{\sigma_r(\lambda)/C^r_n}{\sigma_s(\lambda)/C^s_n}\big]^{\frac{1}{r-s}} .\]

\end{lemma}

\begin{lemma}(See \cite{LRW})\label{lemma3}
(1)Let $u$ be a $k+1$-convex fucntion, $\sigma_k(D^2u)\le C$,
$\lambda=(\lambda_1,\cdots,\lambda_n)$ be the eigenvalues of $D^2u$
with
$\lambda_1\ge\cdots\ge\lambda_n$.
Then there exists a positive constant $K_0$ such that
$\lambda_i+ K_0>0$.

(2)Assume there exists a positive constant $K_0$ such that
$D^2u+ K_0I>0$.
 Let
$\lambda=(\lambda_1,\cdots,\lambda_n)$ be the eigenvalues of $D^2u$,
$\kappa_i=\lambda_i+K_0>0.$ Then
\[\kappa_j\sigma_k^{jj,ii}+\sigma_k^{jj}\ge\sigma_k^{ii}.\]
\end{lemma}

\textbf{Proof.}

(1)\begin{eqnarray*}
C\ge \sigma_k(D^2 u)\ge \lambda_1\sigma_{k-1}(\lambda|1) +\sigma_k(\lambda|1)\ge \lambda_1\cdots\lambda_k\ge \lambda_k^k.
\end{eqnarray*}
Since $\lambda\in \Gamma^k$, we have
\[\lambda_n\ge-\sum_{i=2}^{n-1}\lambda_i \ge-C.\]

(2)
\begin{eqnarray*}
&&\kappa_j\sigma_k^{jj,ii}+\sigma_k^{jj}\\
&=&(K_0+\lambda_j)\sigma_{k-2}(\lambda|ij)+\lambda_i\sigma_{k-2}(\lambda|ij)+\sigma_{k-1}(\lambda|ij)\\
&=&(K_0+\lambda_i)\sigma_{k-2}(\lambda|ij) +\sigma_k^{ii}\ge\sigma_k^{ii}.
\end{eqnarray*}

\begin{lemma}(See \cite{B})\label{lemma4}Let
$B$ be a symmetric matrix, $A$ be a diagonal matrix,
$G$ be a symmetric function of the eigenvalues
of metrices. Let us denote by $\mu(B)$ the
eigenvalues of $B$. Set $g(B)=G(\mu(B))$.
Then
\[
\frac{dg(A+tB)}{ds}\big|_{s=0}
=\frac{\partial G}{\partial\mu_p }B_{pp},\]
and
\[
\frac{d^2g(A+tB)}{ds^2}\big|_{s=0}
=\sum_{p,q}\frac{\partial^2G}{\partial\mu_p\partial\mu_q}B_{pp}B_{qq}
+\sum_{p\neq q}\frac{\frac{\partial G}{\partial\mu_p }-\frac{\partial G}{\partial\mu_q }}{\mu_p-\mu_q}B^2_{pq}.\]
For our case, $G(\mu)=\sum_j\mu_j^m$, $A=D^2u(x_0, t_0)$, $B_{pq}=u_{pqi}$. Then
\[
\frac{dg(A+tB)}{ds}\big|_{s=0}
=m\mu^{m-1}_pu_{ppi},\]
and
\[
\frac{d^2g(A+tB)}{ds^2}\big|_{s=0}
=m(m-1)\sum_p\mu_p^{m-2} u^2_{ppi}
+m\sum_{p\neq q}\frac{\mu_p^{m-1}-\mu_q^{m-1} }{\mu_p-\mu_q}u^2_{pqi}.\]

\end{lemma}

\begin{lemma}(See  \cite{GRW})\label{lemma0}Let~$k>\mu$,~$\alpha=\frac{1}{k-\mu}$.
For any $\delta>0$,
we have
\begin{equation}
-\sigma_k^{pp,qq}u_{pph}u_{qqh}
+(1-\alpha+\frac{\alpha}{\delta})\frac{(\sigma_k)^2_h}{\sigma_k}\ge
\sigma_k(\alpha+1-\delta\alpha)
\big[\frac{(\sigma_\mu)_h}{\sigma_\mu}\big]^2
-\frac{\sigma_k}{\sigma_\mu}\sigma_\mu^{pp,qq}u_{pph}u_{qqh}.\label{2.6789}
\end{equation}

\end{lemma}

\begin{lemma}\label{lemma6}(See \cite{C13}, \cite{HPX})
Suppose that $W$ is diagonal, $W_{11}\geq\cdot\cdot\cdot\geq W_{nn}$
and the eigenvalues of $W$ lie in $\Gamma^2$.
If $\xi_{ij}$ is symmetric and
$$\sum_{i=2}^{n}\sigma_{2}^{ii}\xi_{ii}+\sigma_{2}^{11}\xi_{11}=\eta,$$ then
$$-\sum_{i\neq j}\xi_{ii}\xi_{jj}\geq \frac{n-1}{2\sigma_{2}(W)}
\frac{[2\sigma_{2}(W)\xi_{11}-W_{11}\eta]^2}{[(n-1)W_{11}^2+2(n-2)\sigma_{2}(W)]}
-\frac{\eta^2}{2\sigma_{2}(W)}.$$

For our case, let
$W=D^2u$, $\sigma_2 =-\frac{1}{u_t}$,
$\xi_{ij}=u_{ijl}$, $\eta=\frac{u_{tl}}{u_t^2}$.
Assume $\epsilon$ is a small positive constant and
$u_{11}\ge\sqrt{\frac{2(n-2)}{(n-1)\epsilon m_1}}$.
Then
\begin{eqnarray}
 &&-\sum_{i\neq j}u_{iil}u_{jjl}\notag\\
 &\ge&\frac{ 1}{2\sigma_{2}}
\frac{(2\sigma_{2}u_{11l})^2\big(1-\frac{1}{3(1+\epsilon)+1}\big)+(u_{11}\eta)^2(1-3(1+\epsilon)-1)}{ (1+\epsilon)u_{11}^2 }
-\frac{u_{tl}^2}{2\sigma_{2}u_t^4}\notag\\
 &\geq&\big(\frac{2(1-\frac{1}{3(1+\epsilon)+1})}{1+\epsilon}\big)\sigma_2\frac{u_{11l}^2}{u^2_{11}}
 +\big(\frac{1-(3(1+\epsilon)+1)}{2(1+\epsilon)}-\frac{1}{2}\big) \frac{u_{tl}^2}{u_t^2}\notag\\
 &\geq& \frac{6}{3(1+\epsilon)+1} \sigma_2\frac{u_{11l}^2}{u^2_{11}}
 -2\frac{u_{tl}^2}{u_t^2}.\label{2.99}
\end{eqnarray}

\end{lemma}

\section{A Pogorelov estimate for the \textbf{$k+1$-convex-monotone} solutions to parabolic $k$-Hessian equations }

In this section, we consider Pogorelov estimates for parabolic $k$-Hessian equations \eqref{1.100}. We shall prove Theorem \ref{estimate7}.

Since $u=0$ on $\partial\Omega$, we have $u\leq 0$ in $\Omega$ by the Comparison principle.
By Lemma \ref{lemma3}, there
exists $K_0>0$ such that
$D^2u+ K_0I>0$.
Take the test function
\[\phi=ma^2\log (-u)+ \log P_m+\frac{m}{2}a x_i^2,\]
where $P_m=\sum_i \kappa_i^m$,
$\kappa_i=\lambda_i+K_0>0$.
Constants $a$ and $m$
are positive constants
to be determined later.
Assume the maximum of $\phi$ is attained at $(x_0, t_0)$,
$u_{ij}(x_0, t_0)$ is diagonal and
$u_{11}(x_0, t_0)\ge\cdots\ge u_{nn}(x_0, t_0)$.

Then
\begin{equation}\label{2.1}
0=\frac{1}{m}\phi_i(x_0, t_0)=a^2\frac{u_i}{u}+ \frac{\kappa_l^{m-1}u_{lli}}{P_m}+ ax_i.
\end{equation}

By Lemma \ref{lemma4}, we obtain
\begin{eqnarray}\label{2.2}
0
&\ge& \frac{1}{m}\sigma_k^{ii}\phi_{ii}(x_0, t_0)\notag\\
&=&\sigma_k^{ii}
  \big[a^2\frac{u_{ii}}{u}
  -a^2\frac{u^2_i}{u^2}
  + \frac{\kappa_l^{m-1}u_{llii}}{P_m}
  +\frac{1}{P_m}\sum_{p\neq q}\frac{\kappa^{m-1}_p-\kappa^{m-1}_q}{\kappa_p-\kappa_q}u^2_{pqi}\notag\\
  &&+ \frac{(m-1)\kappa_l^{m-2}u^2_{lli}}{P_m}
  -m\frac{(\kappa_l^{m-1}u_{lli})^2}{P_m^2}
  + a\big]\notag\allowdisplaybreaks\\
&=&\sigma_k^{ii}
  \big[a^2\frac{u_{ii}}{u}
  -a^2\frac{u^2_i}{u^2}
  + \frac{\kappa_l^{m-1}u_{llii}}{P_m}
  +\frac{1}{P_m}\sum_{p\neq q}\sum_{r=0}^{m-2}\kappa_p^r\kappa_q^{m-2-r}u^2_{pqi}\notag\\
  &&+ \frac{(m-1)\kappa_l^{m-2}u^2_{lli}}{P_m}
  -m\frac{(\kappa_l^{m-1}u_{lli})^2}{P_m^2}
  + a\big].\notag\allowdisplaybreaks\\
\end{eqnarray}
Moreover,
\begin{equation}\label{2.3}
0\le \phi_t(x_0, t_0)=a^2\frac{u_t}{u}
+\frac{\kappa_l^{m-1}u_{llt}}{P_m}.
\end{equation}

Now differentiating equations (\ref{1.100}), we obtain
\begin{equation}\label{2.4}
u_{ti}\sigma_k=-u_t  \sigma_k^{jj}u_{jji}
\end{equation}
and

\begin{equation}\label{2.5}
 u_{tii}\sigma_k+ 2u_{ti} \sigma_k^{jj}u_{jji}
=(-u_t)
 \big(\sigma_k^{pq,rs}u_{pqi}u_{rsi}+ \sigma_k^{jj}u_{jjii}\big).
\end{equation}

Note that
\[\sigma_k^{pq,rs}u_{pql}u_{rsl}=\sum_{p\neq q}\sigma_k^{pp,qq}(u_{ppl}u_{qql}-u_{pql}^2).\]
Then (\ref{2.5}) implies that
\begin{equation}\label{2.8}
u_{tll}\sigma_k+ 2u_{tl} \sigma_k^{jj}u_{jjl}
+(-u_t)\sum_{p\neq q}\sigma_k^{pp,qq}(-u_{ppl}u_{qql}+u_{pql}^2)
= (-u_t)\sigma_k^{jj}u_{jjll}.\end{equation}

Then by (\ref{2.1})-(\ref{2.8}),
we have
\begin{eqnarray*}
0
&\ge&
  \sigma_k^{ii}
  \big[a^2\frac{u_{ii}}{u}
  -a^2\frac{u^2_i}{u^2}
  + \frac{\kappa_l^{m-1}u_{llii}}{P_m}
  +\frac{1}{P_m}\sum_{p\neq q}\sum_{r=0}^{m-2}\kappa_p^r\kappa_q^{m-2-r}u^2_{pqi}\allowdisplaybreaks\\
  &&+ (m-1)\frac{\kappa_l^{m-2}u^2_{lli}}{P_m}
  -m\frac{(\kappa_l^{m-1}u_{lli})^2}{P_m^2}
  + a\big]\allowdisplaybreaks\\
&\ge&
  \sigma_k^{ii}
  \big[a^2\frac{u_{ii}}{u}
  -a^2\frac{u^2_i}{u^2}
  +\frac{1}{P_m}\sum_{p\neq q}\sum_{r=0}^{m-2}\kappa_p^r\kappa_q^{m-2-r}u^2_{pqi}\allowdisplaybreaks\\
  &&+ (m-1)\frac{\kappa_l^{m-2}u^2_{lli}}{P_m}
  -m\frac{(\kappa_l^{m-1}u_{lli})^2}{P_m^2}
  + a\big]\allowdisplaybreaks\\
  &&+ \frac{\kappa_l^{m-1}}{P_m}
  \big(\frac{u_{tll}\sigma_k+ 2u_{tl} \sigma_k^{jj}u_{jjl}}{-u_t}
+\sum_{p\neq q}\sigma_k^{pp,qq}(-u_{ppl}u_{qql}+u_{pql}^2)\big)\allowdisplaybreaks\\
&\ge&
  (-a^2\frac{u_t}{u})\cdot\frac{\sigma_k}{-u_t}
  + \sigma_k^{ii}
  \big[a^2\frac{u_{ii}}{u}
  - \frac{1}{a^2}\big(\frac{\kappa_l^{m-1}u_{lli}}{P_m}+a x_i\big)^2\allowdisplaybreaks
  +\frac{1}{P_m}\sum_{p\neq q}\sum_{r=0}^{m-2}\kappa_p^r\kappa_q^{m-2-r}u^2_{pqi}\allowdisplaybreaks\\
  &&+ (m-1)\frac{\kappa_l^{m-2}u^2_{lli}}{P_m}
  -m\frac{(\kappa_l^{m-1}u_{lli})^2}{P_m^2}
  + a\big]
  +\frac{\kappa_l^{m-1}}{P_m}
  \big(\frac{2u_{tl}\sigma_ku_{tl}}{u^2_t}\big)\allowdisplaybreaks\\
   && +\frac{\kappa_l^{m-1}}{P_m}\sum_{p\neq q}\sigma_k^{pp,qq}
  \big(
  -u_{ppl}u_{qql}+u_{pql}^2\big)\allowdisplaybreaks\\
&\ge&
  a^2\frac{\sigma_k}{u}
  + \sigma_k^{ii}
  \big[a^2\frac{u_{ii}}{u}
  - \frac{1}{a^2}\big(\frac{\kappa_l^{m-1}u_{lli}}{P_m}+a x_i\big)^2\allowdisplaybreaks
  +\frac{1}{P_m}\sum_{p\neq q}\sum_{r=0}^{m-2}\kappa_p^r\kappa_q^{m-2-r}u^2_{pqi}\allowdisplaybreaks\\
  &&+ (m-1)\frac{\kappa_l^{m-2}u^2_{lli}}{P_m}
  -m\frac{(\kappa_l^{m-1}u_{lli})^2}{P_m^2}
  + a\big]
  +\frac{\kappa_l^{m-1}}{P_m}
  \big(\frac{ 2u_{tl} \sigma_ku_{tl}}{u^2_t}\big)\allowdisplaybreaks\\
   && +\frac{\kappa_l^{m-1}}{P_m}\sum_{p\neq q}\sigma_k^{pp,qq}
  \big(  -u_{ppl}u_{qql}+ u^2_{pql}\big)\allowdisplaybreaks\\
&\ge&
   a^2\frac{\sigma_k}{u}
  +\sigma_k^{ii}
  \big[a^2\frac{u_{ii}}{u}-2x_i^2 +a
  - \frac{2}{a^2}\big(\frac{\kappa_l^{m-1}u_{lli}}{P_m}\big)^2\allowdisplaybreaks
  +\frac{1}{P_m}\sum_{p\neq q}\sum_{r=0}^{m-2}\kappa_p^r\kappa_q^{m-2-r}u^2_{pqi}\allowdisplaybreaks\\
  &&+ (m-1)\frac{\kappa_l^{m-2}u^2_{lli}}{P_m}
  -m\frac{(\kappa_l^{m-1}u_{lli})^2}{P_m^2}
  \big]+\frac{\kappa_l^{m-1}}{P_m}
  \big(\frac{ 2u_{tl} \sigma_ku_{tl}}{u^2_t}\big)\allowdisplaybreaks\\
   && +\frac{\kappa_l^{m-1}}{P_m}
  \big(  -\sum_{p\neq q}\sigma_k^{pp,qq}u_{ppl}u_{qql}\big)
  +2\sum_{i\neq j}\sigma_k^{ii,jj}\frac{\kappa_j^{m-1}}{P_m} u^2_{ijj}\allowdisplaybreaks\\
&:=&
  Ca^2\frac{\sigma_k}{u}
  + \sigma_k^{ii}
  \big[
  -{2}{}( x_i)^2
  + a\big]+\frac{\kappa_l^{m-1}}{P_m}
  \big(\frac{ 2u_{tl} \sigma_ku_{tl}}{u^2_t}\big)\allowdisplaybreaks\\
   && +\frac{\kappa_l^{m-1}}{P_m}
  \big(  -\sum_{p\neq q}\sigma_k^{pp,qq}u_{ppl}u_{qql} \big)\allowdisplaybreaks
  +\sum_iI_i,\\
\end{eqnarray*}
where
\begin{eqnarray}
I_i&=&\allowdisplaybreaks
\frac{2}{P_m}\sum_{j\neq i}\sum_{r=0}^{m-2}\sigma_k^{jj}\kappa_i^r\kappa_j^{m-2-r}u^2_{ijj}\allowdisplaybreaks
+\sigma_k^{ii} (m-1)\frac{\kappa_j^{m-2}u^2_{jji}}{P_m}\notag\\
&&
-\sigma_k^{ii}(m+\frac{2}{a^2})\frac{(\sum_j\kappa_j^{m-1}u_{jji})^2}{P_m^2}
+2\sum_{i\neq j}\sigma_k^{ii,jj}\frac{\kappa_j^{m-1}}{P_m} u^2_{ijj}.
\end{eqnarray}

We claim that:
\begin{claim}\label{claim}
Suppose $u$ is the $k+1$-convex solution of
(\ref{1.100}) with $m_1\le -u_t\le m_2$. Then,
either
\begin{equation}
|u_{11}|\le C,\label{2356}
\end{equation}
or
\begin{equation}
\frac{\kappa_l^{m-1}}{P_m}
  \big(\frac{ 2u_{tl} \sigma_ku_{tl}}{u^2_t}\big)\allowdisplaybreaks
     +\frac{\kappa_l^{m-1}}{P_m}
  \sigma_k^{pp,qq}\big(  -u_{ppl}u_{qql} \big)\allowdisplaybreaks
  +\sum_iI_i\ge 0.\label{3.100}
  \end{equation}
  \end{claim}
It is obviously that
(\ref{2356}) implies (\ref{1.5}).
If \eqref{3.100} holds,
combining  Lemma \ref{lemma2},
we can obtain
\begin{eqnarray*}
0
&\ge&
  a^2\frac{\sigma_k}{u}+ \sigma_k^{ii}
  \big[-C+a\big]
  \ge \sigma_k^{ii}+\frac{C}{u}\notag\\
&\ge& C\sigma_1^{\frac{1}{k-1}}\sigma_k^{\frac{k-2}{k-1}}+\frac{C}{u}.\allowdisplaybreaks\\
\end{eqnarray*}
Then (\ref{1.5}) is still holds and we completes the proof of Theorem \ref{estimate7}.

Thus the proofs for (\ref{2356}) and  (\ref{3.100}) in Claim \ref{claim}
are  the remaining questions. By Lemma \ref{lemma3}, we obtain
\begin{equation*}2\kappa_j\sigma_k^{jj,ii}+2\sigma_k^{jj}\ge2\sigma_k^{ii}\end{equation*}
and
\begin{equation}
2\kappa_j\sigma_k^{jj,ii}\kappa_j^{m-2}u^2_{jji}
+2\sigma_k^{jj}\kappa_j^{m-2}u_{jji}^2
\ge2\sigma_k^{ii}\kappa_j^{m-2}u^2_{jji}.\label{3.97}\end{equation}
Moreover, by Cauthy inequality, we have
\begin{equation*}2\sum_{p\neq q\neq i}\kappa_p^{m-1}\kappa_q^{m-1}u_{ppi}u_{qqi}
\le 2\sum_{p\neq q\neq i}\kappa_q^{m-2}\kappa_p^m u_{qqi}^2. \end{equation*}
It yields
\begin{equation}
-(\sum_j\kappa_j^{m-1}u_{jji})^2
\ge-\sum_j \kappa_j^{2m-2}u_{jji} ^2
-\sum_{p\neq q\neq i}\kappa_q^{m-2}\kappa_p^m u_{qqi}^2
-2\sum_{p\neq i}\kappa_p^{m-1}\kappa_i^{m-1}u_{ppi}u_{iii}.\label{3.98}
\end{equation}
Therefore, by (\ref{3.97}) and (\ref{3.98}), we have
\begin{eqnarray}
&&P_m^2I_i\notag\allowdisplaybreaks\\
&\ge& \allowdisplaybreaks
  \sigma_k^{ii}\big[
  {2}{P_m}\sum_{j\neq i} \kappa_j^{m-2}u^2_{ijj}\allowdisplaybreaks
  +(m-1)  {\kappa_j^{m-2}u^2_{jji}}{P_m}
  -(m+\frac{2}{a^2}) (\sum_j\kappa_j^{m-1}u_{jji})^2\big]\notag\allowdisplaybreaks\\
  &&+
  {2}{P_m}\sum_{j\neq i}\sum_{r=1}^{m-2}\sigma_k^{jj}\kappa_i^r\kappa_j^{m-2-r}u^2_{ijj}\notag\allowdisplaybreaks\\
&\ge&\sigma_k^{ii}\sum_{j\neq i}
  \big[ \allowdisplaybreaks
  (m+1)  {\kappa_j^{m-2}}{P_m}
  -(m+\frac{2}{a^2}) {\kappa_j^{2m-2}}
  \big]u_{jji}^2 \allowdisplaybreaks\notag\\
  &&
  +{2}{P_m}\sum_{j\neq i}\sum_{r=1}^{m-2}\sigma_k^{jj}\kappa_i^r\kappa_j^{m-2-r}u^2_{ijj}
  \allowdisplaybreaks\notag\\
  &&
  + (m-1)\sigma_k^{ii}{(\kappa_i)^{m-2}}{P_m}u_{iii}^2
  -(m+\frac{2}{a^2})\sigma_k^{ii}{(\kappa_i)^{2m-2}}u_{iii}^2\allowdisplaybreaks\notag\\
  &&-(m+\frac{2}{a^2})\sum_{j\neq l\neq i}\sigma_k^{ii}{(\kappa_j)^{ m-2}\kappa_l^mu^2_{jji}}
  -2(m+\frac{2}{a^2})\sum_{j\neq i}\sigma_k^{ii}{(\kappa_j)^{ m-1}u_{jji}}{(\kappa_i)^{ m-1}u_{iii}}\notag\\
&\ge&\sigma_k^{ii}\sum_{j\neq i}
  \big[ \allowdisplaybreaks
  (m+1)  {\kappa_j^{m-2}}{P_m}
  -(m+\frac{2}{a^2})\sum_{l\neq i} {\kappa_j^{m-2}}\kappa_l^m
  \big]u_{jji}^2\allowdisplaybreaks\notag\\
  &&+
  {2}{P_m}\sum_{j\neq i}\sum_{r=1}^{m-2}\sigma_k^{jj}\kappa_i^r\kappa_j^{m-2-r}u^2_{ijj}
  \allowdisplaybreaks\notag\\
  &&
  + (m-1)\sigma_k^{ii}{(\kappa_i)^{m-2}}{P_m}u_{iii}^2
  -(m+\frac{2}{a^2})\sigma_k^{ii}{(\kappa_i)^{2m-2}}u_{iii}^2\allowdisplaybreaks\notag\\
  &&
  -2(m+\frac{2}{a^2})\sum_{j\neq i}\sigma_k^{ii}{(\kappa_j)^{ m-1}u_{jji}}{(\kappa_i)^{ m-1}u_{iii}}\notag\\
&\ge&
  \sigma_k^{ii}\sum_{j\neq i}
  \big[ (m+1)  {\kappa_j^{m-2}}\kappa_i^m
  +(1-\frac{2}{a^2})\sum_{l\neq i} {(\kappa_j)^{ m-2}\kappa_l^m }
  \big]u_{jji}^2\allowdisplaybreaks\notag\\
  &&+
  {2}{P_m}\sum_{j\neq i}\sum_{r=1}^{m-2}\sigma_k^{jj}\kappa_i^r\kappa_j^{m-2-r}u^2_{ijj}
  \allowdisplaybreaks\notag\\
  &&
  + (m-1)\sigma_k^{ii}{(\kappa_i)^{m-2}}{\kappa_l^m}u_{iii}^2
  -(m+\frac{2}{a^2})\sigma_k^{ii}{(\kappa_i)^{2m-2}}u_{iii}^2\allowdisplaybreaks\notag\\
  &&
  -2(m+\frac{2}{a^2})\sum_{j\neq i}\sigma_k^{ii}{(\kappa_j)^{ m-1}u_{jji}}{(\kappa_i)^{ m-1}u_{iii}}.\label{formuA}
\end{eqnarray}

We divide the proof into two cases: $i>1$ and $i=1$.

(A) $i>1$. In this case, we assert that
\begin{equation}
I_i\ge 0\textrm{ for }i>1.
\label{3.456}
\end{equation}
We further divide case (A) into three subcases to prove the above assertion \eqref{3.456}.

(A1) $\lambda_i\le \lambda_j$, $\lambda_i\ge {K_0}{}$.
\begin{eqnarray*}&&2\sigma_k^{jj}\kappa_i^{m-2-l}\kappa_j^{l}\allowdisplaybreaks\\
&=&2(\lambda_i+\sigma_{k-1}(\lambda|ij))\kappa_i^{m-2-l}\kappa_j^{l}\allowdisplaybreaks\\
&\ge&(\kappa_i+\sigma_{k-1}(\lambda|ij))\kappa_i^{m-2-l}\kappa_j^{l}\allowdisplaybreaks\\
&\ge&(\kappa_j+\sigma_{k-1}(\lambda|ij))\kappa_i^{m-1-l}\kappa_j^{l-1}\allowdisplaybreaks\\
&=&\kappa_i^{m-1-l} \kappa_j^{l-1}\sigma_k^{ii}.
\end{eqnarray*}
Then combining $P_m\ge \kappa_1^m$, we have
\[
  {2}{P_m}\sum_{j\neq i}\sum_{r=1}^{m-2}\sigma_k^{jj}\kappa_i^{m-2-r}\kappa_j^{r}u^2_{ijj}
  \ge{}\sum_{j\neq i}(m-2)\sigma_k^{ii}  \kappa_i^m\kappa_j^{m-2}u^2_{ijj}.
\]

(A2) $\lambda_i\le \lambda_j$, $\lambda_i< {K_0}{}$.
Let $k\le r\le k+6$.
We may assume $\kappa_1\ge\kappa_i^{r+2}$
and $\kappa_1^{k+1}\sigma_k^{jj}\ge\kappa_1^{k}\lambda_1\sigma_k^{11}
\ge \kappa_1\kappa^{k-1}_1C\sigma_k\ge\sigma_k^{ii}$.
Hence, for $2\le r\le 8$, we have
\[ \kappa_1^m\sigma_k^{jj}\kappa_i^{m-2-r}\kappa_j^r
\ge\kappa_1^{k+1}\sigma_k^{jj}\kappa_i^m\kappa_i^{-r-2}\kappa_1\kappa_1^{m-k-2}\kappa_j^r
\ge\sigma_k^{ii}\kappa_j^{m-2}\kappa_i^{m},\]
where we have used $\kappa_1
\ge\kappa_i^{r+2}$.
Thus
\[\sum_{j\neq i}\sigma_k^{jj}
  {P_m}\sum_{r=1}^{m-2}\kappa_i^{m-2-r}\kappa_j^{r}u^2_{ijj}
  \ge7\sum_{j\neq i}\sigma_k^{ii} \kappa_i^m\kappa_j^{m-2}u^2_{ijj}.
\]

(A3) $\lambda_i\ge \lambda_j$.
\[\sigma_k^{jj}\kappa_i^{m-2-r}\kappa_j^r
\ge\kappa_i^{m-2-r}\kappa_j^{r }\sigma_k^{ii}.\]
Combining $P_m\ge \kappa_1^m$, we have
\[
  {2}{P_m}\sum_{j\neq i}\sum_{r=1}^{m-2}\sigma_k^{jj}\kappa_i^r\kappa_j^{m-2-r}u^2_{ijj}
  \ge2\sum_{j\neq i}\sigma_k^{ii}(m-2)  \kappa_i^m\kappa_j^{m-2}u^2_{ijj}.
\]

We choose $m\ge9$. From the above three subcases,
we obtain
\begin{eqnarray*}
&&P_m^2I_i\allowdisplaybreaks\\
&\ge&
  \sigma_k^{ii}
  \Big[
  \sum_{j\neq i}
  \big[ (m+1)  {\kappa_j^{m-2}}\kappa_i^m
  +(1-\frac{2}{a^2})\sum_{l\neq i} {(\kappa_j)^{ m-2}\kappa_l^m }
  \big]u_{jji}^2
  \allowdisplaybreaks\\
  &&+\sum_{j\neq i}
  7\kappa_i^m\kappa_j^{m-2}u^2_{ijj}
  + (m-1) {(\kappa_i)^{m-2}}{\kappa_l^m}u_{iii}^2
  -(m+\frac{2}{a^2}) {(\kappa_i)^{2m-2}}u_{iii}^2\allowdisplaybreaks\\
  &&
  -2(m+\frac{2}{a^2})\sum_{j\neq i} {(\kappa_j)^{ m-1}u_{jji}}{(\kappa_i)^{ m-1}u_{iii}}\Big]\\
&\ge&
  \sigma_k^{ii}
  \Big[
  \sum_{j\neq i}
  \big[ (m+1)  {\kappa_j^{m-2}}\kappa_i^m
  +(1-\frac{2}{a^2})\sum_{l\neq i} {(\kappa_j)^{ m-2}\kappa_l^m }
  \big]u_{jji}^2
  \allowdisplaybreaks\\
  &&+\sum_{j\neq i}
  7\kappa_i^m\kappa_j^{m-2}u^2_{ijj}
  + (m-1)\sum_{j\neq i} {(\kappa_i)^{m-2}}{\kappa_j^m}u_{iii}^2
  -(1+\frac{2}{a^2}) {(\kappa_i)^{2m-2}}u_{iii}^2\allowdisplaybreaks\\
  &&
  -2(m+\frac{2}{a^2})\sum_{j\neq i} {(\kappa_j)^{ m-1}u_{jji}}{(\kappa_i)^{ m-1}u_{iii}}\Big]\\
&\ge&
  \sigma_k^{ii}
  \sum_{j\neq i}\Big[
  (m+8)  {\kappa_j^{m-2}}\kappa_i^m
  u_{jji}^2
  \allowdisplaybreaks\allowdisplaybreaks\\
  &&+ (m-2-\frac{2}{a^2}) {(\kappa_i)^{m-2}}{\kappa_j^m}u_{iii}^2\allowdisplaybreaks
  -2(m+\frac{2}{a^2})  {(\kappa_j)^{ m-1}u_{jji}}{(\kappa_i)^{ m-1}u_{iii}}
  \Big]\allowdisplaybreaks\\
&\ge&0,
\end{eqnarray*}
where we choose
\[(m+8)(m-2-\frac{2}{a^2})\ge(m+\frac{2}{a^2})^2 \]
when $a>1$ and $m$ is sufficiently large.

(B) $i=1$.
In this case, we shall prove that either
(\ref{2356}) holds or
\begin{equation}\label{2.999}
{\kappa_h^{m-1}}{P_m}
  \big(\frac{2u_{th} \sigma_ku_{th}}{u^2_t}\big)\allowdisplaybreaks
     +{\kappa_h^{m-1}}{P_m}
  \sigma_k^{pp,qq}\big(  -u_{pph}u_{qqh} \big) +P_m^2I_1\ge0.
\end{equation}
Now set $\delta>0$ such that $\frac{1}{k-\mu+1}\le \delta<1$.
Here $\mu$ is an integer less than $k$.
Then
\[\max\{1-\frac{1}{k}+\frac{1}{k\delta}, 1+(\frac{1}{\delta}-1)\frac{1}{k-\mu}
\}\le 2\]
Note that, by (\ref{2.6789}) in Lemma \ref{lemma0}, we have
\begin{eqnarray}\label{3.78}
 &&\sigma_k^{pp,qq}\big(  -u_{pph}u_{qqh} \big)
 +\frac{ 2u_{th} \sigma_ku_{th}}{u^2_t}\allowdisplaybreaks
  \notag\allowdisplaybreaks\\
  &\ge&-\sigma_k^{pp,qq}u_{pph}u_{qqh}
+(1-\frac{1}{k}+\frac{1}{k\delta})\frac{(\sigma_k)^2_h}{\sigma_k}\ge 0,\textrm{ for }h>1.
\end{eqnarray}
Besides
\begin{eqnarray}
&&
  \frac{ 2u_{t1} \sigma_ku_{t1}}{u^2_t}\allowdisplaybreaks
     +
  \sigma_k^{pp,qq}\big(  -u_{pp1}u_{qq1} \big)\notag\allowdisplaybreaks\\
&\ge&\sigma_k\Big(1+(1-\delta)\frac{1}{k-\mu}\Big)\big[\frac{(\sigma_\mu)_1}{\sigma_\mu}\big]^2
-\frac{\sigma_k}{\sigma_\mu}\sigma_\mu^{pp,qq}u_{pp1}u_{qq1}.\label{3.09097}
\end{eqnarray}
Therefore, from (\ref{3.78}) and (\ref{3.09097}) we have
\begin{eqnarray}\label{3.0909}
&&\frac{\kappa_h^{m-1}}{P_m}
  \big(\frac{2u_{th} \sigma_ku_{th}}{u^2_t}\big)\allowdisplaybreaks
     +\frac{\kappa_h^{m-1}}{P_m}
  \sigma_k^{pp,qq}\big(  -u_{pph}u_{qqh} \big)\notag\\
   &\ge&\frac{\kappa_1^{m-1}}{P_m}\sigma_k\Big(1+(1-\delta)\frac{1}{k-\mu}\Big)\big[\frac{(\sigma_\mu)_1}{\sigma_\mu}\big]^2
-\frac{\kappa_1^{m-1}}{P_m}\frac{\sigma_k}{\sigma_\mu}\sigma_\mu^{pp,qq}u_{pp1}u_{qq1}.
\end{eqnarray}
 Combining (\ref{formuA}) and (\ref{3.0909}),
by direct calculation, the left hand side of (\ref{2.999}) becomes
\begin{eqnarray}
  && {\kappa_h^{m-1}}{P_m}
  \big(\frac{ 2u_{th} \sigma_ku_{th}}{u^2_t}\big)\allowdisplaybreaks
     + {\kappa_h^{m-1}}{P_m}
  \sigma_k^{pp,qq}\big(  -u_{pph}u_{qqh} \big)
  +  {P^2_m}I_1\allowdisplaybreaks\notag\\
&\ge&{\kappa_h^{m-1}}{P_m}
  \big(\frac{ 2u_{th} \sigma_ku_{th}}{u^2_t}\big)\allowdisplaybreaks
     + {\kappa_h^{m-1}}{P_m}
  \sigma_k^{pp,qq}\big(  -u_{pph}u_{qqh} \big)\notag\\
   &&+ \sigma_k^{11}\sum_{j>1}
  \big[ (m+1)  {\kappa_j^{m-2}}\kappa_1^m
  +(1-\frac{2}{a^2})\sum_{l\neq 1} {(\kappa_j)^{ m-2}\kappa_l^m }
  \big]u_{jj1}^2
  \allowdisplaybreaks\notag\\
  &&+ \sum_{j\neq 1}\sigma_k^{jj}
  {2}P_m\sum_{r=1}^{m-2}\kappa_1^r\kappa_j^{m-2-r}u^2_{1jj}
  +  (m-1)\sigma_k^{11}{(\kappa_1)^{m-2}}{\kappa_l^m}u_{111}^2\allowdisplaybreaks\notag\\
  &&
  - (m+\frac{2}{a^2})\sigma_k^{11}{(\kappa_1)^{2m-2}}u_{111}^2\allowdisplaybreaks
 \notag\\
  &&
  -2(m+\frac{2}{a^2})\sigma_k^{11}\sum_{j\neq 1} {(\kappa_j)^{ m-1}u_{jj1}}{(\kappa_1)^{ m-1}u_{111}}
  \notag\\
&\ge&
  {\kappa_h^{m-1}}{P_m}
  \big(\frac{ 2u_{th} \sigma_ku_{th}}{u^2_t}\big)\allowdisplaybreaks
     + {\kappa_h^{m-1}}{P_m}
  \sigma_k^{pp,qq}\big(  -u_{pph}u_{qqh} \big)\notag\\
  &&+ \sigma_k^{11}\sum_{j>1}
  \big[ (m+3)  {\kappa_j^{m-2}}\kappa_1^m
  \big]u_{jj1}^2
  \allowdisplaybreaks\notag\\
  &&+ \sum_{j\neq 1}\sigma_k^{jj}
  {2}P_m \kappa_1^{m-2}u^2_{1jj}
  + \sum_{j\neq 1} (m-1)\sigma_k^{11}{(\kappa_1)^{m-2}}{\kappa_j^m}u_{111}^2\allowdisplaybreaks\notag\\
  &&
  + \big(-1-\frac{2}{a^2}\big)\sigma_k^{11}{(\kappa_1)^{2m-2}}u_{111}^2\allowdisplaybreaks \notag\\
  &&
  -2(m+\frac{2}{a^2})\sum_{j\neq 1}\sigma_k^{11} {(\kappa_j)^{ m-1}u_{jj1}}{(\kappa_1)^{ m-1}u_{111}}\allowdisplaybreaks\notag\\
&\ge&
  {\kappa_h^{m-1}}{P_m}
  \big(\frac{ 2u_{th} \sigma_ku_{th}}{u^2_t}\big)\allowdisplaybreaks
     + {\kappa_h^{m-1}}{P_m}
  \sigma_k^{pp,qq}\big(  -u_{pph}u_{qqh} \big)\notag\\
  && +\sum_{j\neq 1}\sigma_k^{jj}
  {2}P_m\kappa_1^{m-2}u^2_{1jj}
  + \big(-1-\frac{2}{a^2}\big)\sigma_k^{11}{(\kappa_1)^{2m-2}}u_{111}^2\allowdisplaybreaks \notag\\
  &&+\sigma_k^{11}\sum_{j>1}
  \Big[
  (m+3) {\kappa_j^{m-2}}\kappa_1^m
  u_{jj1}^2
  \allowdisplaybreaks
  +    (m-1) {(\kappa_1)^{m-2}}{\kappa_j^m}u_{111}^2\allowdisplaybreaks\notag\\
  &&-2(m+\frac{2}{a^2})  {(\kappa_j)^{ m-1}u_{jj1}}{(\kappa_1)^{ m-1}u_{111}}
\Big] \allowdisplaybreaks\notag\\
&\ge&{\kappa_h^{m-1}}{P_m}
  \big(\frac{ 2u_{th} \sigma_ku_{th}}{u^2_t}\big)\allowdisplaybreaks
     + {\kappa_h^{m-1}}{P_m}
  \sigma_k^{pp,qq}\big(  -u_{pph}u_{qqh} \big)\notag\\
    && +\sum_{j\neq 1}\sigma_k^{jj}
  {2}P_m\kappa_1^{m-2}u^2_{1jj}
  + \big(-1-\frac{2}{a^2}\big)\sigma_k^{11}{(\kappa_1)^{2m-2}}u_{111}^2\allowdisplaybreaks\notag\\
&\ge&\kappa_1^{m-1}P_m\sigma_k(1+\alpha(1-\delta))\big[\frac{(\sigma_\mu)_1}{\sigma_\mu}\big]^2
   -\kappa_1^{m-1}P_m\frac{\sigma_k}{\sigma_\mu}\sigma_\mu^{pp,qq}u_{pp1}u_{qq1}\notag\\
    && +\sigma_k^{jj}
  {2}P_m \sum_{j\neq 1}\kappa_1^{m-2}u^2_{1jj}
  + \big(-1-\frac{2}{a^2}\big)\sigma_k^{11}{(\kappa_1)^{2m-2}}u_{111}^2\allowdisplaybreaks\notag\\
&:=&E_1+ E_2+ E_3,\label{3.015}
\end{eqnarray}
where we have used
$(m+3)(m-1)\ge(m+\frac{2}{a^2})^2$ for $a>\sqrt{2}$ and $m$ sufficiently large. Here
$$E_1=\kappa_1^{m-1}P_m\sigma_k(1+\alpha(1-\delta))\big[\frac{(\sigma_\mu)_1}{\sigma_\mu}\big]^2
 -\kappa_1^{m-1}P_m\frac{\sigma_k}{\sigma_\mu}\sigma_\mu^{pp,qq}u_{pp1}u_{qq1},$$
$$E_2=\sigma_k^{jj}
  {2}P_m \sum_{j\neq 1}\kappa_1^{m-2}u^2_{1jj}$$
and
$$E_3=
   \big(-1-\frac{2}{a^2}\big)\sigma_k^{11}{(\kappa_1)^{2m-2}}u_{111}^2.$$

Suppose that there
exists $\delta_k>0$ such that
$\frac{\lambda_j}{\lambda_1}\ge \delta_k\textrm{ for all }2\le j\le k$.
Then
\begin{equation}\label{7777}\frac{1}{m_1}\ge \sigma_k\ge\lambda_1\cdots\lambda_k\ge C\lambda_1^k, \end{equation}
and (\ref{2356}) holds.
Therefore, we can
focus on the otherwise situation and
further divide case (B) into two subcases.
For convenience, let us fix $a$
satisfying $a^2\ge\frac{8}{\alpha(1-\delta)}$.

(B1) Fix $\mu=1$ and $\delta_1=1$. Then we can find $\delta_2>0$ such that
(\ref{2.999}) holds when
\begin{eqnarray}
\frac{\lambda_{2}}{\lambda_1}\label{forC}
\le \delta_{2}.
\end{eqnarray}
In fact, by direct calculation, we have
\begin{eqnarray}
E_1&=&\kappa_1^{m-1}P_m\sigma_k(1+\alpha(1-\delta))\big[\frac{(\sigma_1)_1}{\sigma_1}\big]^2
    \notag\allowdisplaybreaks\\
&\ge&\frac{P_m\kappa_1^{m-1}\sigma_k}{\sigma_1^2}\big[(1+\alpha(1-\delta))
     \sum_a(\sigma_1^{aa}u_{aa1})^2\notag\allowdisplaybreaks\\
     &&+ (1+\alpha(1-\delta))\sum_{a\neq b}\sigma_1^{aa}\sigma_{1}^{bb}u_{aa1}u_{bb1}
    \big]\notag\allowdisplaybreaks\\
&\ge&\frac{P_m\kappa_1^{m-1}\sigma_k}{\sigma_1^2}
     \big[(1 +\frac{ \alpha(1-\delta)}{2})\sum_{a\le1}(\sigma_1^{aa}u_{aa1})^2
     -C\sum_{b>1}(\sigma_1^{bb}u_{bb1})^2\big]\notag\allowdisplaybreaks\\
&\ge& { \kappa_1^{2m-2}\sigma^{11}_k}
     \big[(1 +\frac{ \alpha(1-\delta)}{2})(1+\delta_1^m) (1-\frac{C\lambda_{2}}{\lambda_1})^2( u_{111})^2\big]\notag\allowdisplaybreaks\\
     &&-C\frac{P_m\kappa_1^{m-3}\lambda_1^2 }{\sigma_1^2}{}\sum_{b>1}(u_{bb1})^2 \notag\allowdisplaybreaks\\
&\ge& { \kappa_1^{2m-2}\sigma^{11}_k}
     \big[ (1+ \frac{ \alpha(1-\delta)}{4}) ( u_{111})^2\big]
     -C{P_m\kappa_1^{m-3}} \sum_{b>1}( u_{bb1})^2. \label{23.01}
\end{eqnarray}
if $\delta_{2}$ is small enough.

Note that
for $k\ge j>1$,
\begin{equation}
\kappa_1\sigma_{k-1}(\lambda|j)\ge\frac{\lambda_1\cdots\lambda_k\kappa_1}{C\lambda_j}\ge\frac{\sigma_k}{C\delta_{2}}\ge C\frac{1}{\delta_{2}}.
\label{23.6798}
\end{equation}

For $j\ge k+1$,

\begin{equation}
\kappa_1\sigma_{k-1}(\lambda|j)\ge\frac{\lambda_1\lambda_1\cdots\lambda_k}{C\lambda_k}\ge\frac{\sigma_k}{C\delta_{2}}\ge C\frac{1}{\delta_{2}}.
\label{23.6788}
\end{equation}

Then combining (\ref{23.01}), (\ref{23.6798}) and (\ref{23.6788}), (\ref{3.015}) becomes

\begin{eqnarray}\label{23.016}
&&E_1+ E_2+E_3\notag\allowdisplaybreaks\\
&\ge& \Big[{ \kappa_1^{2m-2}\sigma^{11}_k}
     \big[ (1+ \frac{ \alpha(1-\delta)}{4})  ( u_{111})^2\big]\notag\allowdisplaybreaks\\
     &&+\big(-1-\frac{2}{a^2}\big)\sigma_k^{11}{(\kappa_1)^{2m-2}}u_{111}^2
     \Big]\notag\allowdisplaybreaks\\
     &&+P_m\kappa_1^{m-3}\Big[\sum_{j\neq 1}\kappa_1 \sigma_k^{jj}
     {2}u^2_{1jj}-C \sum_{b>1}( u_{bb1})^2
     \Big]\notag\allowdisplaybreaks\\
&\ge&0,
\end{eqnarray}
if we choose  $\delta_{2}$ sufficiently small.
Then we have proved (\ref{2.999}) when (\ref{forC}) holds.

(B2)
Now we assert that we can further find constants $\delta_3,\cdots, \delta_{k}$,
such that (\ref{2.999}) holds when
\begin{eqnarray}
\frac{\lambda_{\mu+1}}{\lambda_1}\label{3B1}
\le \delta_{\mu+1},
\end{eqnarray}
and
\begin{eqnarray}\frac{\lambda_\mu}{\lambda_1}
\ge \delta_{\mu}\label{3B2}\end{eqnarray}
for some $\mu\in\{2, \cdots, k-1\}$.

To this end, we will prove it by induction.
In other words, we assume (\ref{3B2}) holds firstly.
Then we shall
find $\delta_{\mu+1}>0$ sufficiently small such that
(\ref{2.999}) holds provided we have (\ref{3B1}).
Since $\lambda\in\Gamma_{k+1}\subset\Gamma_{\mu+2}$, we have, for $a, b\le\mu$,
\begin{equation}\label{3.a}
\sigma_\mu^{aa}\ge\frac{\lambda_1\cdots\lambda_\mu}{\lambda_a},
\end{equation}
\begin{equation}\label{3.b}
\sigma_{\mu-1}(\lambda|ab)\le C\frac{\lambda_1\cdots\lambda_{\mu+1}}{\lambda_a\lambda_b},
\end{equation}
\begin{equation}\label{3.c}
\sigma_{\mu}(\lambda|ab)\le C\frac{\lambda_1\cdots\lambda_{\mu+2}}{\lambda_a\lambda_b},
\end{equation}
\begin{equation}\label{3.d}
\sigma_{\mu-2}(\lambda|ab)\le C\frac{\lambda_1\cdots\lambda_{\mu}}{\lambda_a\lambda_b}.
\end{equation}
Combining (\ref{3.a})-(\ref{3.d}), by direct calculation,
we have
\begin{equation}\label{3.e}
0<\sigma_\mu^{aa}\sigma_\mu^{bb}-\sigma_\mu\sigma_\mu^{aa,bb}
=\sigma^2_{\mu-1}(\lambda|ab)-\sigma_{\mu}(\lambda|ab)\sigma_{\mu-2}(\lambda|ab)
\le C(\frac{\lambda_{\mu+1}}{\lambda_b}\sigma_\mu^{aa})^2.
\end{equation}
Moreover, by (\ref{3.e}), we obtain
\begin{eqnarray}
&&\sum_{a\neq b, a,b\le \mu}(\sigma_\mu^{aa}\sigma_\mu^{bb}-\sigma_\mu\sigma_\mu^{aa,bb})u_{aa1}u_{bb1}\notag\allowdisplaybreaks\\
&\ge&-C\sum_{a\neq b, a,b\le \mu}\big(\frac{\lambda_{\mu+1}}{\lambda_b}\big)^2(\sigma_\mu^{aa})^2\notag\allowdisplaybreaks\\
&\ge&-C\sum_{a\le \mu}\big(\frac{\lambda_{\mu+1}}{\delta_\mu\lambda_1}\big)^2(\sigma_\mu^{aa})^2\notag\allowdisplaybreaks\\
&\ge&-C\frac{\delta^2_{\mu+1}}{\delta^2_\mu}\sum_{a\le \mu} (\sigma_\mu^{aa})^2 \notag\allowdisplaybreaks\\
&\ge&-\frac{ 1}{a^2}\sum_{a\le \mu} (\sigma_\mu^{aa})^2 \allowdisplaybreaks, \label{3.21}
\end{eqnarray}
if $\delta_{\mu+1}$ is sufficiently small.
Besides,
\begin{eqnarray}
&&2\sum_{a\neq b, a \le \mu, b>\mu}(\sigma_\mu^{aa}\sigma_\mu^{bb}-\sigma_\mu\sigma_\mu^{aa,bb})u_{aa1}u_{bb1}\notag\allowdisplaybreaks\\
&\ge&-\frac{ 1}{a^2}\sum_{a\le \mu}(\sigma_\mu^{aa}u_{aa1})^2- {a^2}\sum_{b>\mu}(\sigma_\mu^{bb}u_{bb1})^2,\label{3.22}
\end{eqnarray}
and
\begin{eqnarray}
&&\sum_{a\neq b, a , b>\mu}(\sigma_\mu^{aa}\sigma_\mu^{bb}-\sigma_\mu\sigma_\mu^{aa,bb})u_{aa1}u_{bb1}\notag\allowdisplaybreaks\\
&\ge&- \sum_{b>\mu}(\sigma_\mu^{bb}u_{bb1})^2\label{3.23}
\end{eqnarray}

From (\ref{3.21})-(\ref{3.23}), if we
choose $a$ sufficiently large
and $\delta_{\mu+1}$ sufficiently small, $E_1$ becomes

\begin{eqnarray}
E_1&=&\kappa_1^{m-1}P_m\sigma_k(1+\alpha(1-\delta))\big[\frac{(\sigma_\mu)_1}{\sigma_\mu}\big]^2
   -\kappa_1^{m-1}P_m\frac{\sigma_k}{\sigma_\mu}\sigma_\mu^{pp,qq}u_{pp1}u_{qq1}\notag\allowdisplaybreaks\\
&\ge&\frac{P_m\kappa_1^{m-1}\sigma_k}{\sigma_\mu^2}\big[(1+\alpha(1-\delta))
     \sum_a(\sigma_\mu^{aa}u_{aa1})^2
     + \alpha(1-\delta)\sum_{a\neq b}\sigma_\mu^{aa}\sigma_{\mu}^{bb}u_{aa1}u_{bb1}\notag\allowdisplaybreaks\\
    &&+\sum_{a\neq b}(\sigma_\mu^{aa}\sigma_\mu^{bb}-\sigma_\mu\sigma_\mu^{aa,bb})u_{aa1}u_{bb1}\big]\notag\allowdisplaybreaks\\
&\ge&\frac{P_m\kappa_1^{m-1}\sigma_k}{\sigma_\mu^2}
     \big[(1 +\frac{ \alpha(1-\delta)}{2})\sum_{a\le\mu}(\sigma_\mu^{aa}u_{aa1})^2
     -C\sum_{b>\mu}(\sigma_\mu^{bb}u_{bb1})^2\big]\notag\allowdisplaybreaks\\
&\ge& { \kappa_1^{2m-2}\sigma^{11}_k}
     \big[(1 +\frac{ \alpha(1-\delta)}{2})(1+\delta_\mu^m)\sum_{a\le\mu}(1-\frac{C\lambda_{\mu+1}}{\lambda_a})^2( u_{aa1})^2\big]\notag\allowdisplaybreaks\\
     &&-C\frac{P_m\kappa_1^{m-3}\lambda_1^2 }{\sigma_\mu^2}{}\sum_{b>\mu}(\sigma_\mu^{bb}u_{bb1})^2 \notag\allowdisplaybreaks\\
&\ge& { \kappa_1^{2m-2}\sigma^{11}_k}
     \big[(1 +\frac{ \alpha(1-\delta)}{2})(1+\delta_\mu^m)\sum_{a\le\mu}(1-\frac{C\lambda_{\mu+1}}{\delta_\mu\lambda_1})^2( u_{aa1})^2\big]\notag\allowdisplaybreaks\\
     &&-C\frac{P_m\kappa_1^{m-3}\lambda_1^2}{\lambda_\mu^2}\sum_{b>\mu}( u_{bb1})^2 \notag\allowdisplaybreaks\\
&\ge& { \kappa_1^{2m-2}\sigma^{11}_k}
     \big[ (1+ \frac{ \alpha(1-\delta)}{4})\sum_{a\le\mu} ( u_{aa1})^2\big]
     -C\frac{P_m\kappa_1^{m-3}}{\delta_\mu^2}\sum_{b>\mu}( u_{bb1})^2, \label{3.01}
\end{eqnarray}
if $\delta_{\mu+1}$ is sufficiently small.

Note that
for $k\ge j>\mu$,
\begin{equation}
\kappa_1\sigma_{k-1}(\lambda|j)\ge\frac{\lambda_1\cdots\lambda_k\kappa_1}{C\lambda_j}\ge\frac{\sigma_k}{C\delta_{\mu+1}}\ge C\frac{1}{\delta_{\mu+1}}.
\label{3.6798}
\end{equation}

For $j\ge k+1$,

\begin{equation}
\kappa_1\sigma_{k-1}(\lambda|j)\ge\frac{\lambda_1\lambda_1\cdots\lambda_k}{C\lambda_k}\ge\frac{\sigma_k}{C\delta_{\mu+1}}\ge C\frac{1}{\delta_{\mu+1}}.
\label{3.6788}
\end{equation}

Then combining (\ref{3.01})-(\ref{3.6788}) and (\ref{3.015}) we have

\begin{eqnarray}\label{3.016}
&&E_1+ E_2+E_3\notag\allowdisplaybreaks\\
&\ge& \Big[{ \kappa_1^{2m-2}\sigma^{11}_k}
     \big[ (1+ \frac{ \alpha(1-\delta)}{4})\sum_{a\le\mu} ( u_{aa1})^2\big]\notag\\
     &&+\big(-1-\frac{2}{a^2}\big)\sigma_k^{11}{(\kappa_1)^{2m-2}}u_{111}^2
     \Big]\notag\allowdisplaybreaks\\
     &&+P_m\kappa_1^{m-3}\Big[\sum_{j\neq 1}\kappa_1 \sigma_k^{jj}
     {2}u^2_{1jj}-C\frac{1}{\delta_\mu^2}\sum_{b>\mu}( u_{bb1})^2
     \Big]\notag\allowdisplaybreaks\\
&\ge&0,
\end{eqnarray}
if we choose  $\delta_{\mu+1}$ sufficiently small.
Hence, combining (\ref{7777}), (\ref{3.015}), (\ref{23.016}) and (\ref{3.016}), we have proved that
either (\ref{2356}) or (\ref{2.999}) holds.
Besides,
from (A1)-(A3), we have (\ref{3.456}).
Thus we have proved the Claim
\ref{claim} and
we complete the proof of Theorem \ref{estimate7}.

\section{A Pogorelov estimate for the 3-convex-monotone solutions to
parabolic 2-Hessian equations}
In this section, we shall prove the Theorem \ref{estimate1}.
Since $u=0$ on $\partial\Omega$, we have $u\leq 0$ in $\Omega$ by the Comparison principle (see Theorem 17.1 in Page 443 of  \cite{GT}).
By Lemma \ref{lemma3}, there
exists $K_0>0$ such that
$D^2u+ K_0I>0$.
Take the test function
\[\phi=ma^2\log (-u)+ \log P_m+\frac{m}{2}a x_i^2+ \frac{mN}{2}|Du|^2,\]
where $P_m=\sum_i \kappa_i^m$,
$\kappa_i=\lambda_i+K_0>0$.
Constants $a$ and $m$
are positive constants
to be determined later.
Assume the maximum of $\phi$ is attained at $(x_0, t_0)$,
$u_{ij}(x_0, t_0)$ is diagonal and
$u_{11}(x_0, t_0)\ge\cdots\ge u_{nn}(x_0, t_0)$.

Then
\begin{equation}\label{82.1}
0=\frac{1}{m}\phi_i(x_0, t_0)=a^2\frac{u_i}{u}+ \frac{\kappa_l^{m-1}u_{lli}}{P_m}+ ax_i +N u_k u_{ki}.
\end{equation}

By Lemma \ref{lemma4} we obtain
\begin{eqnarray}\label{82.2}
0
&\ge& \frac{1}{m}\sigma_2^{ii}\phi_{ii}(x_0, t_0)\notag\\
&=&\sigma_2^{ii}
  \big[a^2\frac{u_{ii}}{u}
  -a^2\frac{u^2_i}{u^2}
  + \frac{\kappa_l^{m-1}u_{llii}}{P_m}
  +\frac{1}{P_m}\sum_{p\neq q}\frac{\kappa^{m-1}_p-\kappa^{m-1}_q}{\kappa_p-\kappa_q}u^2_{pqi}\notag\\
  &&+ \frac{(m-1)\kappa_l^{m-2}u^2_{lli}}{P_m}
  -m\frac{(\kappa_l^{m-1}u_{lli})^2}{P_m^2}
  + a+ Nu_{ki}^2 + Nu_ku_{kii}\big]\notag\allowdisplaybreaks\\
&=&\sigma_2^{ii}
  \big[a^2\frac{u_{ii}}{u}
  -a^2\frac{u^2_i}{u^2}
  + \frac{\kappa_l^{m-1}u_{llii}}{P_m}
  +\frac{1}{P_m}\sum_{p\neq q}\sum_{r=0}^{m-2}\kappa_p^r\kappa_q^{m-2-r}u^2_{pqi}\notag\\
  &&+ \frac{(m-1)\kappa_l^{m-2}u^2_{lli}}{P_m}
  -m\frac{(\kappa_l^{m-1}u_{lli})^2}{P_m^2}
  + a+ Nu_{ki}^2 + Nu_ku_{kii}\big].\notag\allowdisplaybreaks
\end{eqnarray}
Moreover,
\begin{equation}\label{82.3}
0\le \phi_t(x_0, t_0)=a^2\frac{u_t}{u}
+\frac{\kappa_l^{m-1}u_{llt}}{P_m} + Nu_ku_{kt}.
\end{equation}

Now differentiating equation (\ref{1.100}), we obtain
\begin{equation}\label{82.4}
u_{ti}\sigma_2=-u_t  \sigma_2^{jj}u_{jji}
\end{equation}
and

\begin{equation}\label{82.5}
 u_{tii}\sigma_2+ 2u_{ti} \sigma_2^{jj}u_{jji}
=(-u_t)
 \big(\sigma_2^{pq,rs}u_{pqi}u_{rsi}+ \sigma_2^{jj}u_{jjii}\big).
\end{equation}

Note that
\[\sigma_2^{pq,rs}u_{pql}u_{rsl}=u_{ppl}u_{qql}-u_{pql}^2.\]
Then (\ref{82.5}) implies that
\begin{equation}\label{82.8}
u_{tll}\sigma_2+ 2u_{tl} \sigma_2^{jj}u_{jjl}
+(-u_t)(-u_{ppl}u_{qql}+u_{pql}^2)
= (-u_t)\sigma_2^{jj}u_{jjll}.\end{equation}

Then by (\ref{82.1})-(\ref{82.8})
and (\ref{2.99}), 
we have

\begin{eqnarray*}
0
&\ge&
  \sigma_2^{ii}
  \big[a^2\frac{u_{ii}}{u}
  -a^2\frac{u^2_i}{u^2}
  + \frac{\kappa_l^{m-1}u_{llii}}{P_m}
  +\frac{1}{P_m}\sum_{p\neq q}\sum_{r=0}^{m-2}\kappa_p^r\kappa_q^{m-2-r}u^2_{pqi}\allowdisplaybreaks\\
  &&+ (m-1)\frac{\kappa_l^{m-2}u^2_{lli}}{P_m}
  -m\frac{(\kappa_l^{m-1}u_{lli})^2}{P_m^2}
  + a+ Nu_{ki}^2 + Nu_ku_{kii}\big]\allowdisplaybreaks\\
&\ge&
  \sigma_2^{ii}
  \big[a^2\frac{u_{ii}}{u}
  -a^2\frac{u^2_i}{u^2}
  +\frac{1}{P_m}\sum_{p\neq q}\sum_{r=0}^{m-2}\kappa_p^r\kappa_q^{m-2-r}u^2_{pqi}\allowdisplaybreaks\\
  &&+ (m-1)\frac{\kappa_l^{m-2}u^2_{lli}}{P_m}
  -m\frac{(\kappa_l^{m-1}u_{lli})^2}{P_m^2}
  + a+ Nu_{ki}^2 \big]+ Nu_k\big(\frac{2u_{tk}\sigma_2}{-u_t}\big)\allowdisplaybreaks\\
  &&+ \frac{\kappa_l^{m-1}}{P_m}
  \big(\frac{u_{tll}\sigma_2+ 2u_{tl} \sigma_2^{jj}u_{jjl}}{-u_t}
+(-u_{ppl}u_{qql}+u_{pql}^2)\big)\allowdisplaybreaks\\
&\ge&
  (-a^2\frac{u_t}{u})\cdot\frac{\sigma_2}{-u_t}
  + \sigma_2^{ii}
  \big[a^2\frac{u_{ii}}{u}
  - \frac{1}{a^2}\big(\frac{\kappa_l^{m-1}u_{lli}}{P_m}+a x_i\big)^2\allowdisplaybreaks
  +\frac{1}{P_m}\sum_{p\neq q}\sum_{r=0}^{m-2}\kappa_p^r\kappa_q^{m-2-r}u^2_{pqi}\allowdisplaybreaks\\
  &&+ (m-1)\frac{\kappa_l^{m-2}u^2_{lli}}{P_m}
  -m\frac{(\kappa_l^{m-1}u_{lli})^2}{P_m^2}
  + a\big]
  +\frac{\kappa_l^{m-1}}{P_m}
  \big(\frac{ 2u_{tl}\sigma_2u_{tl}}{u^2_t}\big)\allowdisplaybreaks\\
   && +\frac{\kappa_l^{m-1}}{P_m}
  \big(
  -u_{ppl}u_{qql}+u_{pql}^2\big)+N\sigma_2^{ii}u_{ii}^2\allowdisplaybreaks\\
&\ge&
  a^2\frac{\sigma_2}{u}
  + \sigma_2^{ii}
  \big[a^2\frac{u_{ii}}{u}
  - \frac{1}{a^2}\big(\frac{\kappa_l^{m-1}u_{lli}}{P_m}+a x_i\big)^2\allowdisplaybreaks
  +\frac{1}{P_m}\sum_{p\neq q}\sum_{r=0}^{m-2}\kappa_p^r\kappa_q^{m-2-r}u^2_{pqi}\allowdisplaybreaks\\
  &&+ (m-1)\frac{\kappa_l^{m-2}u^2_{lli}}{P_m}
  -m\frac{(\kappa_l^{m-1}u_{lli})^2}{P_m^2}
  + a\big]
  +\frac{\kappa_l^{m-1}}{P_m}
  \big(\frac{ 2u_{tl} \sigma_2u_{tl}}{u^2_t}\big)\allowdisplaybreaks\\
   && +\frac{\kappa_l^{m-1}}{P_m}
  \big( \sigma_2 \frac{6}{3(1+\epsilon)+1} \frac{u_{11l}^2}{u^2_{11}}-2\frac{u_{tl}^2}{u_t^4}+ u^2_{pql}\big)
  +N\sigma_2^{ii}u_{ii}^2\allowdisplaybreaks\\
&\ge&
   a^2\frac{\sigma_2}{u}
  + \sigma_2^{ii}
  \big[a^2\frac{u_{ii}}{u}
  - \frac{1}{a^2}\big(\frac{\kappa_l^{m-1}u_{lli}}{P_m}+a x_i\big)^2\allowdisplaybreaks
  +\frac{1}{P_m}\sum_{p\neq q}\sum_{r=0}^{m-2}\kappa_p^r\kappa_q^{m-2-r}u^2_{pqi}\allowdisplaybreaks\\
  &&+ (m-1)\frac{\kappa_l^{m-2}u^2_{lli}}{P_m}
  -m\frac{(\kappa_l^{m-1}u_{lli})^2}{P_m^2}
  + a\big]
  \allowdisplaybreaks\\
  && +\frac{\kappa_l^{m-1}}{P_m}
  \big(
  \sigma_2 \frac{6}{3(1+\epsilon)+1} \frac{u_{11l}^2}{u^2_{11}}+u_{pql}^2 \big)
  +N\sigma_2^{ii}u_{ii}^2\allowdisplaybreaks\\
&\ge&
   a^2\frac{\sigma_2}{u}
  +\sigma_2^{ii}
  \big[a^2\frac{u_{ii}}{u}-2x_i^2 +a
  - \frac{2}{a^2}\big(\frac{\kappa_l^{m-1}u_{lli}}{P_m}\big)^2\allowdisplaybreaks
  +\frac{1}{P_m}\sum_{p\neq q}\sum_{r=0}^{m-2}\kappa_p^r\kappa_q^{m-2-r}u^2_{pqi}\allowdisplaybreaks\\
  &&+ (m-1)\frac{\kappa_l^{m-2}u^2_{lli}}{P_m}
  -m\frac{(\kappa_l^{m-1}u_{lli})^2}{P_m^2}
  \big]+\frac{\kappa_l^{m-1}}{P_m}
  \big(
  \sigma_2\frac{6}{3(1+\epsilon)+1} \frac{u_{11l}^2}{u^2_{11}}+u_{pql}^2 \big)
  +N\sigma_2^{ii}u_{ii}^2
  \allowdisplaybreaks\\
&\ge&
  Ca^2\frac{\sigma_2}{u}
  + \sigma_2^{ii}
  \big[
  -{2}{}( x_i)^2
  + a\big]+
  \frac{\kappa_1^{m-1}}{P_m}\sigma_2 \frac{6}{3(1+\epsilon)+1}   \frac{u_{111}^2}{u^2_{11}}
  +N\sigma_2^{ii}u_{ii}^2\allowdisplaybreaks
  +\sum_iI_i,\\
\end{eqnarray*}
where
\begin{eqnarray}
I_i&=&\allowdisplaybreaks
\frac{2}{P_m}\sum_{j\neq i}\sum_{r=0}^{m-2}\sigma_2^{jj}\kappa_i^r\kappa_j^{m-2-r}u^2_{ijj}\allowdisplaybreaks
+\sigma_2^{ii} (m-1)\frac{\kappa_j^{m-2}u^2_{jji}}{P_m}\notag\\
&&
-\sigma_2^{ii}(m+\frac{2}{a^2})\frac{(\sum_j\kappa_j^{m-1}u_{jji})^2}{P_m^2}
+\frac{\kappa_l^{m-1}}{P_m}
\big( u_{pql}^2 \big).
\end{eqnarray}

Now we assert that
\begin{equation}
\sum_l\frac{\kappa_1^{m-1}}{P_m}
\big[ \sigma_2 \frac{6}{3(1+\epsilon)+1}   \frac{u_{111}^2}{u^2_{11}}
\big]+ \sum_iI_i\ge 0.\label{83.100}
  \end{equation}
If (\ref{83.100}) holds,
we obtain
\begin{eqnarray}
0
&\ge&
  a^2\frac{\sigma_2}{u}+ \sigma_2^{ii}
  \big[-C+a+ Nu_{ii}^2\big]
  .\allowdisplaybreaks\label{4.999}
\end{eqnarray}
 Note that $\sigma_2^{11}u_{11}\ge c_0$.
 Thus (\ref{4.999}) becomes
 \begin{eqnarray}
 \frac{C}{-u}+ Cu_{11}(C-a)
  &\ge& c_0Nu_{11},
  \allowdisplaybreaks\label{4.998}
\end{eqnarray}
where $N$ is sufficiently large.
Then \[-uu_{11}\le C.\]
It completes the proof of Theorem \ref{estimate1}.

Now the remaining question is the proof of assertion (\ref{83.100}). Note that
by Lemma \ref{lemma3}, we obtain
\begin{equation*}2\kappa_j+2\sigma_2^{jj}\ge2\sigma_2^{ii}\end{equation*}
and
\begin{equation}
2\kappa_j\kappa_j^{m-2}u^2_{jji}
+2\sigma_2^{jj}\kappa_j^{m-2}u_{jji}^2
\ge2\sigma_2^{ii}\kappa_j^{m-2}u^2_{jji}\label{83.97}\end{equation}
Moreover,
 by Cauthy inequality, we have
\begin{equation*}2\sum_{p\neq q\neq i}\kappa_p^{m-1}\kappa_q^{m-1}u_{ppi}u_{qqi}
\le 2\sum_{p\neq q\neq i}\kappa_q^{m-2}\kappa_p^m u_{qqi}^2. \end{equation*}
It yields
\begin{equation}
-(\sum_j\kappa_j^{m-1}u_{jji})^2
\ge-\sum_j \kappa_j^{2m-2}u_{jji} ^2
-\sum_{p\neq q\neq i}\kappa_q^{m-2}\kappa_p^m u_{qqi}^2
-2\sum_{p\neq i}\kappa_p^{m-1}\kappa_i^{m-1}u_{ppi}u_{iii}.\label{83.98}
\end{equation}
Therefore, by (\ref{83.97}) and (\ref{83.98}) we have
\begin{eqnarray*}
&&P_m^2I_i\allowdisplaybreaks\\
&\ge& \allowdisplaybreaks
  \sigma_2^{ii}\big[
  {2}{P_m}\sum_{j\neq i} \kappa_j^{m-2}u^2_{ijj}\allowdisplaybreaks
  +(m-1)  {\kappa_j^{m-2}u^2_{jji}}{P_m}
  -(m+\frac{2}{a^2}) (\sum_j\kappa_j^{m-1}u_{jji})^2\big]\allowdisplaybreaks\\
  &&+
  {2}{P_m}\sum_{j\neq i}\sum_{r=1}^{m-2}\sigma_2^{jj}\kappa_i^r\kappa_j^{m-2-r}u^2_{ijj}\allowdisplaybreaks\\
&\ge&\sigma_2^{ii}\sum_{j\neq i}
  \big[ \allowdisplaybreaks
  (m+1)  {\kappa_j^{m-2}}{P_m}
  -(m+\frac{2}{a^2}) {\kappa_j^{2m-2}}
  \big]u_{jji}^2
  \allowdisplaybreaks\\
  &&+
  {2}{P_m}\sum_{j\neq i}\sum_{r=1}^{m-2}\sigma_2^{jj}\kappa_i^r\kappa_j^{m-2-r}u^2_{ijj}
  + (m-1)\sigma_2^{ii}{(\kappa_i)^{m-2}}{P_m}u_{iii}^2
  -(m+\frac{2}{a^2})\sigma_2^{ii}{(\kappa_i)^{2m-2}}u_{iii}^2\allowdisplaybreaks\\
  &&-(m+\frac{2}{a^2})\sum_{j\neq l\neq i}\sigma_2^{ii}{(\kappa_j)^{ m-2}\kappa_l^mu^2_{jji}}
  -2(m+\frac{2}{a^2})\sum_{j\neq i}\sigma_2^{ii}{(\kappa_j)^{ m-1}u_{jji}}{(\kappa_i)^{ m-1}u_{iii}}\\
&\ge&\sigma_2^{ii}\sum_{j\neq i}
  \big[ \allowdisplaybreaks
  (m+1)  {\kappa_j^{m-2}}{P_m}
  -(m+\frac{2}{a^2})\sum_{l\neq i} {\kappa_j^{m-2}}\kappa_l^m
  \big]u_{jji}^2
  \allowdisplaybreaks\\
  &&+
  {2}{P_m}\sum_{j\neq i}\sum_{r=1}^{m-2}\sigma_2^{jj}\kappa_i^r\kappa_j^{m-2-r}u^2_{ijj}
  + (m-1)\sigma_2^{ii}{(\kappa_i)^{m-2}}{P_m}u_{iii}^2
  -(m+\frac{2}{a^2})\sigma_2^{ii}{(\kappa_i)^{2m-2}}u_{iii}^2\allowdisplaybreaks\\
  &&
  -2(m+\frac{2}{a^2})\sum_{j\neq i}\sigma_2^{ii}{(\kappa_j)^{ m-1}u_{jji}}{(\kappa_i)^{ m-1}u_{iii}}\\
&\ge&
  \sigma_2^{ii}\sum_{j\neq i}
  \big[ (m+1)  {\kappa_j^{m-2}}\kappa_i^m
  +(1-\frac{2}{a^2})\sum_{l\neq i} {(\kappa_j)^{ m-2}\kappa_l^m }
  \big]u_{jji}^2
  \allowdisplaybreaks\\
  &&+
  {2}{P_m}\sum_{j\neq i}\sum_{r=1}^{m-2}\sigma_2^{jj}\kappa_i^r\kappa_j^{m-2-r}u^2_{ijj}
  + (m-1)\sigma_2^{ii}{(\kappa_i)^{m-2}}{\kappa_l^m}u_{iii}^2
  -(m+\frac{2}{a^2})\sigma_2^{ii}{(\kappa_i)^{2m-2}}u_{iii}^2\allowdisplaybreaks\\
  &&
  -2(m+\frac{2}{a^2})\sum_{j\neq i}\sigma_2^{ii}{(\kappa_j)^{ m-1}u_{jji}}{(\kappa_i)^{ m-1}u_{iii}}.
\end{eqnarray*}

We divide the proof into two different cases: $i>1$ and $i=1$.

(A)$i>1$. In this case, we shall prove
\begin{equation}
I_i>0, \textrm{ for  } i>1.\label{789}
\end{equation}

We further divide case (A) into three subcases.

(A1) $\lambda_i\le \lambda_j$, $\lambda_i\ge {K_0}{}$.
\begin{eqnarray*}&&2\sigma_2^{jj}\kappa_i^{m-2-l}\kappa_j^{l}\allowdisplaybreaks\\
&\ge&(\kappa_j+\sigma_1(\lambda|ij))\kappa_i^{m-1-l}\kappa_j^{l-1}\allowdisplaybreaks\\
&=&\kappa_i^{m-1-l} \kappa_j^{l-1}\sigma_2^{ii}.
\end{eqnarray*}
Then combining $P_m\ge \kappa_1^m$, we have
\[
  {2}{P_m}\sum_{j\neq i}\sum_{r=1}^{m-2}\sigma_2^{jj}\kappa_i^{m-2-r}\kappa_j^{r}u^2_{ijj}
  \ge{}\sum_{j\neq i}(m-2)\sigma_2^{ii}  \kappa_i^m\kappa_j^{m-2}u^2_{ijj}.
\]

(A2)$\lambda_i\le \lambda_j$, $\lambda_i< {K_0}{}$.
Let $2\le r\le 8$.
We may assume $\kappa_1\ge\kappa_i^{r+2}$
and $\kappa_1^{3}\sigma_2^{jj}\ge\kappa_1^{2}\lambda_1\sigma_2^{11}
\ge \kappa_1\kappa_1C\sigma_2\ge\sigma_2^{ii}$.
Hence, for $2\le r\le 8$, we have
\[ \kappa_1^m\sigma_2^{jj}\kappa_i^{m-2-r}\kappa_j^r
\ge\kappa_1^{3}\sigma_2^{jj}\kappa_i^m\kappa_i^{-r-2}\kappa_1\kappa_1^{m-4}\kappa_j^r
\ge\sigma_2^{ii}\kappa_j^{m-2}\kappa_i^{m}.\]
Thus
\[\sum_{j\neq i}\sigma_2^{jj}
  {P_m}\sum_{r=1}^{m-2}\kappa_i^{m-2-r}\kappa_j^{r}u^2_{ijj}
  \ge7\sum_{j\neq i}\sigma_2^{ii} \kappa_i^m\kappa_j^{m-2}u^2_{ijj}.
\]

(A3)$\lambda_i\ge \lambda_j$.
\[\sigma_2^{jj}\kappa_i^{m-2-r}\kappa_j^r
\ge\kappa_i^{m-2-r}\kappa_j^{r }\sigma_k^{ii}.\]
Combining $P_m\ge \kappa_1^m$, we have
\[
  {2}{P_m}\sum_{j\neq i}\sum_{r=1}^{m-2}\sigma_2^{jj}\kappa_i^r\kappa_j^{m-2-r}u^2_{ijj}
  \ge2\sum_{j\neq i}\sigma_2^{ii}(m-2)  \kappa_i^m\kappa_j^{m-2}u^2_{ijj}.
\]

From the above three subcases,
we have
\begin{eqnarray*}
&&P_m^2I_i\allowdisplaybreaks\\
&\ge&
  \sigma_2^{ii}
  \sum_{j\neq i}\Big[
  (m+8)  {\kappa_j^{m-2}}\kappa_i^m
  u_{11i}^2
  \allowdisplaybreaks\allowdisplaybreaks\\
  &&+ (m-2-\frac{2}{a^2}) {(\kappa_i)^{m-2}}{\kappa_j^m}u_{iii}^2\allowdisplaybreaks
  -2(m+\frac{2}{a^2})  {(\kappa_j)^{ m-1}u_{jj1}}{(\kappa_i)^{ m-1}u_{iii}}
  \Big]\allowdisplaybreaks\\
&\ge&0,
\end{eqnarray*}
where we have used
\[(m+8)(m-2-\frac{2}{a^2})\ge(m+\frac{2}{a^2})^2 \]
when $m$ is large and $a>1$.

(B)$i=1$. In this case,
we shall prove that
\begin{equation}\label{797}
\frac{(\kappa_1)^{m-1}}{P_m}\sigma_2 \frac{6}{3(1+\epsilon)+1} \frac{u_{111}^2}{u^2_{11}}
  +  I_1\ge0.
\end{equation}
In fact,
it is easy to see that
\begin{eqnarray*}
  &&{(\kappa_1)^{m-1}}{P_m}\sigma_2 \frac{6}{3(1+\epsilon)+1} \frac{u_{111}^2}{u^2_{11}}
=   \frac{6}{3(1+\epsilon)+1}  {(\kappa_1)^{m-1}}{P_m}\sigma_2\frac{u_{111}^2}{\lambda_1^2}\\
&\ge&
   \frac{6}{3(1+\epsilon)+1}  {(\kappa_1)^{m-2}}{P_m}\sigma^{11}_2 {u_{111}^2},
\end{eqnarray*}
where we have used
\begin{equation}\sigma_2\ge\lambda_1\sigma_2^{11}.
\label{3.900}\end{equation}
Note that
\begin{equation}\frac{6}{3(1+\epsilon)+1} >\frac{5}{4}
\label{3.901}\end{equation}
when $\epsilon$ is small. Then
\begin{eqnarray*}
  &&{(\kappa_1)^{m-1}}{P_m}\sigma_2 \frac{6}{3(1+\epsilon)+1} \frac{u_{111}^2}{u^2_{11}}
  +  {P^2_m}I_1\allowdisplaybreaks\\
&\ge& \sigma_2^{11}\sum_{j>1}
  \big[ (m+1)  {\kappa_j^{m-2}}\kappa_1^m
  +(1-\frac{2}{a^2})\sum_{l\neq 1} {(\kappa_j)^{ m-2}\kappa_l^m }
  \big]u_{jj1}^2
  \allowdisplaybreaks\\
  &&+ \sigma_2^{11}
  {2}P_m \sum_{j\neq 1}\sum_{r=0}^{m-3}\kappa_1^r\kappa_j^{m-2-r}u^2_{1jj}
  +  (m-1)\sigma_2^{11}{(\kappa_1)^{m-2}}{\kappa_l^m}u_{111}^2\allowdisplaybreaks\\
  &&
  - (m+\frac{2}{a^2})\sigma_2^{11}{(\kappa_1)^{2m-2}}u_{111}^2\allowdisplaybreaks
  + \frac{6}{3(1+\epsilon)+1}  {(\kappa_1)^{m-2}}{P_m}\sigma^{11}_2 {u_{111}^2}\\
  &&
  -2(m+\frac{2}{a^2})\sum_{j\neq 1} {(\kappa_j)^{ m-1}u_{jj1}}{(\kappa_1)^{ m-1}u_{111}}
  \\
&\ge& \sigma_2^{11}\sum_{j>1}
  \big[ (m+1)  {\kappa_j^{m-2}}\kappa_1^m
  \big]u_{jj1}^2
  \allowdisplaybreaks\\
  &&+ \sigma_2^{11}
  {2}P_m \sum_{j\neq 1}\sum_{r=0}^{m-3}\kappa_1^r\kappa_j^{m-2-r}u^2_{1jj}
  + \sum_{j\neq 1} (m+\frac{1}{4})\sigma_2^{11}{(\kappa_1)^{m-2}}{\kappa_j^m}u_{111}^2\allowdisplaybreaks\\
  &&
  + \big(\frac{5}{4}-1-\frac{2}{a^2}\big)\sigma_2^{11}{(\kappa_1)^{2m-2}}u_{111}^2\allowdisplaybreaks \\
  &&
  -2(m+\frac{2}{a^2})\sum_{j\neq 1}\sigma_2^{11} {(\kappa_j)^{ m-1}u_{jj1}}{(\kappa_1)^{ m-1}u_{111}}\allowdisplaybreaks\\
&\ge&
  \sigma_2^{11}\sum_{j>1}
  \Big[
  (m+1) {\kappa_j^{m-2}}\kappa_1^m
  u_{jj1}^2
  \allowdisplaybreaks
  +    (m+\frac{1}{4}) {(\kappa_1)^{m-2}}{\kappa_j^m}u_{111}^2\allowdisplaybreaks\\
  &&-2(m+\frac{2}{a^2})  {(\kappa_j)^{ m-1}u_{jj1}}{(\kappa_1)^{ m-1}u_{111}}
\Big]\ge0,\allowdisplaybreaks\\
\end{eqnarray*}
if $m$ is sufficiently large and $a^2>\frac{16}{5}$.
From (\ref{789}) and (\ref{797}),
we have proved the assertion (\ref{83.100}).
 $\Box$
\begin{remark}
We point out a fact. The power $8$ in (\ref{100}) can be improved
to any constant larger than $4$. Indeed, for any
$a^2>4$, let $1<\eta<\frac{a^2}{4}$,
$0<\epsilon<\frac{4(\eta-1)}{3+6\eta}$.
Then we have
\[\frac{6}{3(1+\epsilon)+1}-1-\frac{2}{a^2}\ge0.\]
So the argument is still valid.
\end{remark}

\section{Proof of the Liouville theorems}
\textbf{Proof of the Liouville theorem for parabolic k-Hessian equations:}
In this section, we give the proof of Theorem \ref{maintheorem}.
The proof is classical.

Let $u$ be a solution of
equation (\ref{4.10}).
Set
\[v(x, t)=\frac{u(Rx, R^2t)-R^2}{R^2},\]
and
\[\Omega_R=\{(x, t)|u(Rx, R^2t)-R^2\le 0\}.\]
Then $v$ satisfies
\begin{equation}\label{1.1}
\left\{
\begin{array}{rl}
-v_t \sigma_k(D^2v) =1, \ \textrm{in}\  \Omega_R,\\
v=0, \ \textrm{on}\  \partial \Omega_R.
\end{array}
\right.
\end{equation}
Note that
\[b|Rx|^2-c\le u(Rx,0)\le u(Rx, R^2t)\le R^2,\]
and therefore
\[|x|^2\le \frac{1+c}{b}.\]
Thus
$\Omega_R(t)$ is bounded and
\[v(x, t)=\frac{u(Rx, R^2t)-R^2}{R^2}\ge \frac{b|Rx|^2-cR^2-R^2}{R^2}\ge-C.\]
By Theorem \ref{estimate7}, it yields
\[(-v)^\beta(\Delta v)\le C,\]
where $C$ is an absolutely constant.
Besides,
set $\widetilde{\Omega}_{R}=\{(x, t)|u(Rx, R^2t)-\frac{R^2}{2}\le 0\}$.
It is obviously that
$-v\ge\frac{1}{2}$ in $\widetilde{\Omega}_{R}$.
Thus
\[\Delta v\le 2^\beta C.\]
It follows that
\[\Delta u\le 2^\beta C\]
in $\{(y,s)|u(y, s)\le\frac{R^2}{2}\}$. Here
$R$ is arbitrary.
Furthermore,
using Evans-Krylov theory (see \cite{GT}), we obtain
\[\lim_{R\rightarrow\infty}|D^2u|_{C^\alpha}\le\lim_{R\rightarrow\infty} C\frac{|D^2u|_{C^0}}{R^\alpha}
\le\lim_{R\rightarrow\infty}\frac{C}{R^\alpha}=0.\]

It proves the theorem \ref{maintheorem}.
 $\Box$

\end{document}